\documentclass[11pt,a4paper]{article}
\usepackage[ukrainian,english,russian]{babel}
\usepackage{amssymb}
\usepackage{amsmath}
\begin{document}
\title{Some comparison theorems in Finsler-Hadamard manifolds}
\author{A. A. Borisenko, E. A. Olin}
\date{23 May 2006}
\maketitle
\begin{center}
{\it Geometry Department, Math.-Mech. Faculty, Kharkov National
University, Pl. Svoboda 4, 310077-Kharkov, Ukraine.\\}
 E-mail: Alexander.A.Borisenko@univer.kharkov.ua, evolin@mail.ru
\end{center}

\textbf{Abstract}
 We give upper and lower bounds for the ratio of
the volume of metric ball to the area of the metric sphere in
Finsler-Hadamard manifolds with pinched S-curvature. We apply
these estimates to find the limit at the infinity for this ratio.
Derived estimates are the generalization of the well-known result
in Riemannian geometry. We also estimate the volume growth entropy
for the balls in such manifolds.



\section{Introduction}
Finsler geometry is an important generalization of Riemannian
geometry. It was introduced by P. Finsler in 1918 from
considerations of regular problems in the calculus of variations.
In Finsler geometry metric needs not to be quadratic on the
tangent spaces, thus the structure of Finsler spaces is much more
complicated then the structure of Riemannian spaces. But many
notions and theorems were generalized to Finsler geometry from
Riemannian geometry.

In [1], [2] the following result was proved.

\textbf{Theorem 1.}

\textit{Let $M^{n+1}$ be an $(n+1)$-dimensional Hadamard manifold
with the sectional curvature $K$ such that $-k_2^2\leqslant K
\leqslant -k_1^2$, $k_1,k_2>0$. Let $\Omega$ be a compact
$\lambda$-convex domain in $M^{n+1}$ (i.e. domain, whose boundary
is a regular hypersurface with all normal curvatures that are
greater or equal than $\lambda$) with $\lambda\leqslant k_2$. Then
there exist functions $\alpha(r)$ of the inradius and $\beta(R)$
of the circumradius such that $\alpha(r)\rightarrow 1/(nk_1)$ and
$\beta(R)\rightarrow 1/(nk_2)$ as $r$ and $R$ go to infinity and
that
$$\alpha(r)\frac{\lambda}{k_2}\leqslant \frac{Vol(\Omega)}{Vol(\partial
\Omega)}\leqslant\beta(R).$$ }

\textit{As a consequence, for a family $\{\Omega(t)\}_{t \in
\mathbb{R}^+}$ of compact $\lambda$-convex domains with
$\lambda\leqslant k_2$ expanding over the whole space we obtain
$$\frac{\lambda}{nk_2^2}\leqslant \lim_{t\rightarrow\infty}\inf\frac{Vol(\Omega(t))}{Vol(\partial \Omega(t))}\leqslant
\lim_{t\rightarrow\infty}\sup\frac{Vol(\Omega(t))}{Vol(\partial
\Omega(t))}\leqslant \frac{1}{nk_1}.$$}

Our goal is to generalize this theorem for Finsler manifolds. We
consider metric balls as the family $\{\Omega(t)\}_{t \in
\mathbb{R}^+}$. We shall also need bounds for one of
non-Riemannian curvatures, namely S-curvature. As the result we
prove the following theorem.

\textbf{Theorem 2.}

\textit{Let $(M^{n+1},F)$ be an $(n+1)$-dimensional
Finsler-Hadamard manifold that satisfies the following conditions:
\begin{enumerate}
\item Flag curvature satisfies the inequalities $-k_2^2\leqslant K
\leqslant -k_1^2$, $k_1,k_2>0$, \item $S$-curvature satisfies the
inequalities $n\delta_1\leqslant S \leqslant n\delta_2$ such that
$\delta_i<k_i.$
\end{enumerate}
Let $B_r^{n+1}(p)$ be the metric ball of radius $r$ in $M^{n+1}$
with the center at a point $p \in M^{n+1}$, $S_r^{n}(p)=\partial
B_r^{n+1}(p)$ be the metric sphere. Let $Vol=\int dV_F$ be the
measure of Busemann-Hausdorff, $Area=\int dA_F$ the induced
measure on $S_r^{n}(p)$. Then there exist functions $f(r)$ and
$\mathcal{F}(r)$ such that $f(r)\rightarrow 1/(n(k_2-\delta_2))$
and $\mathcal{F}(r)\rightarrow 1/(n(k_1-\delta_1))$ as $r$ goes to
infinity and that
$$f(r) \leqslant \frac{Vol( B_r^{n+1}(p))}{Area(S_r^{n}(p))}\leqslant \mathcal{F}(r).$$
 Here} $$f(r)=\frac{1}{(1-e^{-2k_2r})^n}
\left(\frac{1}{n(k_2-\delta_2)}-\frac{n}{n(k_2-\delta_2)-2k_2}(e^{-2k_2r}-e^{-nr(k_2-\delta_2)})\right)$$
$$\mathcal{F}(r)=\frac{1}{n(k_1-\delta_1)}(1-e^{-nr(k_1-\delta_1)}).$$

\textit{As a consequence, for a family
$\{B_r^{n+1}(p)\}_{r\geqslant0}$ we have
$$\frac{1}{n(k_2-\delta_2)}\leqslant \lim_{r\rightarrow\infty}\inf\frac{Vol( B_r^{n+1}(p))}{Area(S_r^{n}(p))}\leqslant
\lim_{r\rightarrow\infty}\sup\frac{Vol(
B_r^{n+1}(p))}{Area(S_r^{n}(p))}\leqslant
\frac{1}{n(k_1-\delta_1)}.$$} \textit{If $(M^{n+1},F)$ is a space
of constant flag curvature $K=-k^2$ and S-curvature $S=n\delta$,
$\delta<k$, we have}
$$\lim_{r\rightarrow\infty}\frac{Vol(B_r^{n+1}(p))}{Area(S_r^{n}(p))}=\frac{1}{n(k-\delta)}$$

For a Riemannian space $S=0$ and Theorem 2 is a special case of
Theorem 1.

In the section 4 we give estimates for the volume growth entropy
of the balls.
\section{Preliminaries}

In this section we recall some basic facts and theorems from
Finsler geometry that we need. See [3], [4], [5] for details.

\subsection{Finsler Metrics} By definition, a Finsler metric on a manifold is a
family of Minkowski norms on the tangent spaces. A
\textit{Minkowski norm} on a vector space $V^n$ is a nonnegative
function $F:V^n\rightarrow [0,\infty)$ with the following
properties
\begin{enumerate}
\item $F$ is positively homogeneous of degree one, i.e., for any
$y\in V^n$ and any $\lambda>0$, $F(\lambda y)=\lambda F(y)$;

\item $F$ is $C^{\infty}$ on $V^n\backslash\{0\}$ and for any
vector $y\in V^n$ the following  bilinear symmetric form
$g_y:V^n\times V^n\rightarrow \mathbb{R}$ is positive definite,
$$g_y(u,v):=\frac{1}{2}\frac{\partial^2}{\partial t \partial s}
\lbrack F^2(y+su+tv)\rbrack |_{s=t=0}.$$
\end{enumerate}

The property 2. is also called the \textit{strong convexity
property}.

The Minkowski norm is said to be \textit{reversible} if
$F(y)=F(-y)$, $y \in V^n$. In this paper, Minkowski norms are not
assumed to be reversible.

By 1. and 2., one can show that $F(y)>0$ for $y\neq0$ and
$F(u+v)\leqslant F(u)+F(v)$. See [4] for a proof.

A vector space $V^n$ with the Minkowski norm is called a
\textit{Minkowski space}. Notice that reversible Minkowski spaces
are finite-dimensional Banach spaces.

Let $(V^n,F)$ be the Minkowski space. Then the set $I=F^{-1}(1)$
is called the \textit{indicatrix} in the Minkowski space. It is
also called the \textit{unit sphere}.

A set $U \subset V^n$ is said to be \textit{strongly convex} if
there exist a function $F$ satisfying 2. such that $\partial
U=F^{-1}(1)$. Remark that a strong convexity is equivalent to a
positivity of all normal curvatures of $\partial U$ for any
euclidean metric on $V^n$.
\par Let $M^n$ be an $n$-dimensional connected
$C^{\infty}$-manifold. Denote by $TM^n=\bigsqcup_{x\in M^n}T_xM^n$
the tangent bundle of $M^n$, where $T_xM^n$ is the tangent space
at $x$. A \textit{Finsler metric} on $M^n$ is a function
$F:TM^n\rightarrow [0,\infty)$ with the following properties
\begin{enumerate}
\item $F$ is $C^{\infty}$ on $TM^n\backslash\{0\}$;

\item At each point $x \in M^n$, the restriction $F|_{T_xM^n}$ is
a Minkowski norm on $T_xM^n$.
\end{enumerate}

The pair $(M^n,F)$ is called \textit{a Finsler manifold}.
\par
Let $(M^n,F)$ be a Finsler manifold. Let $(x^i,y^i)$ be a standard
local coordinate system in $TM^n$, i.e., $y^i$ are determined by
$y=y^i \frac{\partial}{\partial x^i}|_x$. For a non-zero vector
$y=y^i \frac{\partial}{\partial x^i}$, put
$g_{ij}(x,y):=\frac{1}{2}[F^2]_{y^iy^j}(x,y)$. The induced inner
product $g_y$ is given by $$g_y(u,v)=g_{ij}(x,y)u^iv^j,$$ where
$u=u^i \frac{\partial}{\partial x^i}|_x$,
$v=v^i\frac{\partial}{\partial x^i}|_x$.

By the homogeneity of $F$, we have
$F(x,y)=\sqrt{g_y(y,y)}=\sqrt{g_{ij}(x,y)y^iy^j}$.

In the Riemannian case $g_{ij}$ are functions of $x \in M^n$ only,
and in the Minkowski case $g_{ij}$ are functions of $y \in
T_xM^n=V^n$ only.

\subsection {Measuring of Area}
The notions of length and area are also generalized to Finsler
geometry.

Given a Finsler metric $F$ on a manifold $M^n$.

Let $\{e_i\}_{i=1}^n$ be an arbitrary basis for $T_xM^n$ and
$\{\theta^i\}_{i=1}^n$ the dual basis for $T_x^*M^n$. The set
$$B_x^n=\left\{ (y^i)\in\mathbb{R}^n : F(x, y^ie_i)<1 \right \}$$ is
an open strongly convex open subset in $\mathbb{R}^n$, bounded by
the indicatrix in $T_xM^n$. Then define
$$dV_F=\sigma_F(x)\theta^1\wedge...\wedge\theta^n,$$ where
$$\sigma_F(x):=\frac{Vol_E(\mathbb{B}^n)}{Vol_E(B_x^n)}.$$ Here
$Vol_E(A)$ denotes the Euclidean volume of $A$, and $\mathbb{B}^n$
is the standard unit ball in $\mathbb{R}^n$.

The volume form $dV_F$ determines a regular measure $Vol_F=\int
dV_F$ and is called the \textit{Busemann-Hausdorff volume form}.

For any Riemannian metric $g_{ij}(x)u^iv^j$ the Busemann-Hausdorff
volume form is the standard Riemannian volume form
$$dV_g=\sqrt{\det(g_{ij})}\theta^1\wedge...\wedge\theta^n.$$

Let $\varphi :N^{n-1}\rightarrow M^n$ be a hypersurface in
$(M^n,F)$.

The Finsler metric $F$ determines a local normal vector field as
follows. A vector $n_x$ is called the \textit{normal vector} to
$N^{n-1}$ at $x \in N^{n-1}$ if $n_x \in T_{\varphi(x)}M^n$ and
$g_{n_x}(y,n_x)=0$ for all $y \in T_xN^{n-1}$. It was proved in
[4] that such vector exists. Notice that in general non-symmetric
case the vector $-n_x$ is not a normal vector.

Define now an induced volume form on $N^{n-1}$. Let $n$ be a unit
normal vector field along $N^{n-1}$. Let $\overline{F}=\varphi^*F$
be the induced Finsler metric on $N^{n-1}$ and $dV_{\overline{F}}$
be the Busemann-Hausdorff volume form of $\overline{F}$. For $x
\in N^{n-1}$ we define
$$\zeta(x,n_x):=\frac{Vol_E(\mathbb{B}^n)}{Vol_E(B_x^n)}\frac{Vol_E(B_x^{n-1}(n_x))}{Vol_E(\mathbb{B}^{n-1})}.$$
Here $B_x^n=\left\{ (y^i)\in\mathbb{R}^n : F(y^ie_i)<1 \right \}$.
To define $B_x^{n-1}(n_x)$ we take a basis $\{e_i\}_{i=1}^n$ for
$T_{\varphi(x)}M^n$ such that $e_1=n_x$ and $\{e_i\}_{i=2}^n$ is a
basis for $T_xN^{n-1}$. Then $B_x^{n-1}(n_x)=\left\{
(y^j)\in\mathbb{R}^{n-1} : F(y^je_j)<1\right \}$, where the index
$j$ passes from $2$ to $n$.

Note that if $F$ is Riemannian, then $\zeta \equiv1$.

Set $$dA_F:=\zeta(x,n_x) dV_{\overline{F}}.$$ The form $dA_F$ is
called the \textit{induced volume form} of $dV_F$ with respect to
$n$ [4].

The sense of defining such volume form is given by the
\textit{co-area formula} [4]. We shall need the co-area formula in
one simple case for metric balls:
\begin{equation}
Vol(B(r,p))=\int_0^rVol(S(t,p))dt.
\end{equation}
Here $Vol(S(t,p))$ is the induced volume on $S(t,p)$.

\subsection{Geodesics, Connections and Curvature}
Locally minimizing curves in a Finsler space are determined by a
system of second order differential equations (geodesic
equations).

Let $(M^n,F)$ be a Finsler space, and $c:[a,b]\rightarrow M^n$ a
constant speed piecewise $C^\infty$ curve $F(c,\dot{c})=const$.
Denote the local functions $G^i(x,y)$ by $$G^i(x,y)=\frac{1}{4}
g^{il}(x,y)\left\{2\frac{\partial g_{jl}}{\partial
x^k}(x,y)-\frac{\partial g_{jk}}{\partial x^l}(x,y)
\right\}y^jy^k.$$ We call $G^i(x,y)$ the \textit{geodesic
coefficients} [4]. Notice that in Riemannian case
$G^i(x,y)=\frac{1}{2}\Gamma_{ik}^j(x)y^iy^k$.

Consider the functions $N_j^i(x,y)=\frac{\partial G^i}{\partial
y^j}(x,y)$. They are called the \textit{connection coefficients}.
At each point $x \in M^n$, define a mapping $$D:T_xM^n \times
C^\infty(TM^n)\rightarrow T_xM^n$$ by
$$D_yU:=\{dU^i(y)+U^jN_j^i(x,y)\}\frac{\partial}{\partial x^i}|_x,$$
where $y\in T_xM^n$ and $U\in C^\infty(TM^n)$. We call $D_yU(x)$
the \textit{covariant derivative} of $U$ at $x$ in the direction
$y$.

If $c$ is a solution of the system $D_{\dot{c}}\dot{c}=0$ then it
is called a \textit{geodesic}.

Next, we introduce a notion of curvature in Finsler geometry. At
first, we consider the generalization of Riemann curvature. In
1926, L. Berwald extended the Riemann curvature to Finsler
metrics.

Let $(M^n,F)$ be a Finsler space. For a vector $y\in T_xM^n
\backslash \{0\}$ consider the functions
$$R_i^k(y)=2\frac{\partial G^i}{\partial x^k}-\frac{\partial^2 G^i}{\partial x^j \partial
y^k}y^j+2G^j\frac{\partial^2G^i}{\partial y^j \partial
y^k}-\frac{\partial G^i}{\partial y^j}\frac{\partial G^j}{\partial
y^k}.$$ For every vector $y\in T_xM^n\backslash\{0\}$, define a
linear transformation $$R_y=R_i^k(y)\frac{\partial}{\partial
x^i}\otimes dx^k|_x.$$ Then the family of transformations
$$R=\{R_y:T_xM^n\rightarrow T_xM^n,y\in T_xM^n\backslash\{0\}, x \in
M^n\}$$is called the \textit{Riemann curvature} [4].

Let $P\subset T_xM^n$ be a tangent plane. For a vector $y \in
P\backslash\{0\}$, define
$$K(P,y):=\frac{g_y(R_y(u),u)}{g_y(y,y)g_y(u,u)-g_y(y,u)^2},$$
where $u\in P$ such that $P=span\{y,u\}$. $K(P,y)$ is independent
of $u\in P$. The number $K(P,y)$ is called the \textit{flag
curvature} of the flag $(P,y)$ in $T_xM^n$.

The flag curvature is a generalization of the sectional curvature
in Riemannian geometry. It can be defined in another way. For a
vector $y \in T_xM^n\backslash\{0\}$ consider the Riemannian
metric $\hat{g}(u,v)=g_Y(u,v)$. Here the vector field $Y$ is an
arbitrary extension of the vector $y$. Then the flag curvature
$K(P,y)$ of the flag $(P,y)$ in the Finsler metric $F$ is equal to
the sectional curvature of the plane $P$ in the metric
$\hat{g}(u,v)$ . If we change $y$, then $\hat{g}(u,v)$ and
$K(P,y)$ will also change [3].

Define the \textit{Ricci curvature} by
$$Ric(y)=\sum_{i=1}^nR_i^i(y).$$

A simply-connected Finsler space with non-positive flag curvature
is called a \textit{Finsler-Hadamard space}. In these spaces the
generalization of Cartan-Hadamard's theorem holds [6].

The notions of exponential map, completeness, cut-locus, conjugate
and focal points in Finsler geometry are defined by the same way
as in Riemannian geometry. For details, see [4].

Finally, we introduce some more functions which are called
\textit{non-Riemannian curvatures}. These curvatures all vanish
for Riemannian spaces. We shall need only one of this curvatures,
which is closely connected to the volume form.

Let $(M^n,F)$ be a Finsler space. Consider the Busemann-Hausdorff
volume form $dV_F$ with the density $\sigma_F$. We define
$$\tau(x,y)=\ln\frac{\sqrt{\det(g_{ij}(x,y))}}{\sigma_F(x)},\ y\in T_xM^n.$$
$\tau$ is called the \textit{distortion} of $(M^n,F)$. The
condition $\tau\equiv const$ implies $F$ is a Riemannian metric
[4].

To measure the rate of changes of the distortion along geodesics,
we can define [3,4,5]
$$S(x,y)=\frac{d}{dt}\left[\tau(c(t),\dot{c}(t))\right]|_{t=0}, \ y\in T_xM^n$$ where
$c(t)$ is the geodesic with $\dot{c}(0)=y$. $S$ is called the
\textit{S-curvature}. It is also called the \textit{mean
covariation} and \textit{mean tangent curvature}. The local
formula for the $S$-curvature is
$$S(x,y)=N_m^m(x,y)-\frac{y^m}{\sigma_F(x)} \frac{\partial \sigma_F}{\partial
x^m}(x).$$ One can easily show that $S=0$ for any Riemannian
metric.

A Finsler metric $F$ is said to be of \textit{constant
S-curvature} $\delta$ if $$S(x,y)=\delta F(x,y)$$ for all $y \in
T_xM^n\backslash\{0\}$ and $x \in M^n$. The upper and lower bounds
of S-curvature are defined by the same way.

\subsection{Geometry of Hypersurfaces and Comparison Theorems}
Let $(M^{n},F)$ be a Finsler manifold and
$\varphi:N^{n-1}\rightarrow M^n$ be a hypersurface. Let
$\overline{F}=\varphi^*F$ denote the induced Finsler metric on
$N^{n-1}$. Let $\rho$ be a $C^{\infty}$-distance function on an
open subset $U\subset M^n$ such that $\rho^{-1}(s)=N^{n-1}\cap U$
for some $s$. Let $dV_F$ denotes the Busemann-Hausdorff volume
form of $F$, $dA_t$ denote the induced volume form of
$N_t^{n-1}=\rho^{-1}(t)$. Let $c(t)$ be an integral curve of
$\nabla\rho$ with $c(0)\in N_s^{n-1}$. We have $\rho(c(t))=t$,
hence $c(\varepsilon)\in N_{s+\varepsilon}^{n-1}$ for small
$\varepsilon>0$. By definition, the flow $\phi_{\varepsilon}$ of
$\nabla\rho$ satisfies
$$\phi_{\varepsilon}(c(s))=c(s+\varepsilon).$$
$$\phi_{\varepsilon}:N^{n-1}\cap U=N_s^{n-1}\rightarrow
N_{s+\varepsilon}^{n-1}.$$

The $(n-1)$-form $\phi_{\varepsilon}^*dA_{s+\varepsilon}$ is a
multiply of $dA_s$. Thus there is a function
$\Theta(x,\varepsilon)$ on $N^{n-1}$ such that
$$\phi_{\varepsilon}^*dA_{s+\varepsilon}|_x=\Theta(x,\varepsilon)dA_s|_x,
,\ \forall x \in N^{n-1},$$ $$\Theta(x,0)=1,\ \forall x \in
N^{n-1}.$$ Set $$\Pi_{n_x}=\frac{\partial}{\partial
\varepsilon}\left(\ln \Theta(x,\varepsilon)
\right)|_{\varepsilon=0}.$$ $\Pi_{n_x}$ is called he mean
curvature of $N^{n-1}$ at $x$ with respect to $n_x:=\nabla\rho_x$
[4].

We also need some estimates on the mean curvature of the metric
sphere. The following theorem gives these estimates. For a given
real $\lambda$, put
$$s_{\lambda}(t)=\frac{\sin(\sqrt{\lambda}t)}{\sqrt{\lambda}}, \ \lambda>0 ,$$
$$s_{\lambda}(t)=t, \ \lambda=0,$$
$$s_{\lambda}(t)=\frac{\sinh(\sqrt{-\lambda}t)}{\sqrt{-\lambda}}, \ \lambda<0.$$
\textbf{Theorem 3.} [4]

\textit{Let $(M^n,F)$ be an $n$-dimensional positively complete
Finsler space. Let $\Pi_t$ denote the mean curvature of $S(p,t)$
in the cut-domain of $p$ with respect to the outward-pointing
normal vector.}
\begin{enumerate}
\item \textit{Suppose that $$K\leqslant\lambda,\ S\leqslant
(n-1)\delta.$$ Then \begin{equation}\Pi_t\geqslant
(n-1)\frac{s_{\lambda}'(t)}{s_{\lambda}(t)}-(n-1)\delta.\end{equation}}
\item \textit{Suppose that $$Ric\geqslant n\lambda,\ S\geqslant
-(n-1)\delta.$$ Then \begin{equation}\Pi_t\leqslant
(n-1)\frac{s_{\lambda}'(t)}{s_{\lambda}(t)}+(n-1)\delta.\end{equation}}
\end{enumerate}

\textbf{Theorem 4.} [4] \textit{Let $(M^n,F)$ be an
$n$-dimensional positively complete Finsler space. Suppose that
for constants $\lambda\leqslant0$ and $\delta\geqslant0$ with
$\sqrt{-\lambda}-\delta>0$, the flag curvature and the
$S$-curvature satisfy the inequalities
$$K\leqslant\lambda,\ S\leqslant (n-1)\delta.$$ Then for any regular domain $\Omega\subset M^n,$ $$Vol(\Omega)\leqslant \frac{Vol(\partial \Omega)}{(n-1)(\sqrt{-\lambda}-\delta)}$$}
Remark that the right asymptotic estimate in Theorem 2 is proved
in Theorem 4.

\section{Relation between area and volume for balls in Finsler-Hadamard manifolds}

In this section we prove Theorem 2.

P r o o f $\ $ o f $\ $ T h e o r e m 2 :

Let $S_pM^{n+1}$ denote the unit sphere in $T_pM^{n+1}$. Fix a
vector $y \in S_pM^{n+1}$. Let $\{e_i\}_{i=1}^{n+1}$ be a basis
for $T_pM^{n+1}$ such that
$$e_1=y, \ g_y(y,e_i)=0, \ i=2,...,n+1.$$

Extend $\{e_i\}_{i=1}^n$ to a global frame on $T_pM^{n+1}$ in a
natural way. Let $\{\theta^i\}_{i=1}^{n+1}$ denote the basis for
$T_x^*M^{n+1}$ dual to $\{e_i\}_{i=1}^{n+1}$. Express $dV_F$ at
$p$ by
$$dV_F(p)=\sigma_F(p)\theta^1\wedge...\wedge\theta^{n+1},$$
$$\sigma_F(p)=\frac{Vol_E(\mathbf{B}^{n+1})}{Vol_E(\left\{ (y^i)\in\mathbb{R}^{n+1} : F(y^ie_i)<1 \right
\})}.$$ Thus we obtain the volume form $dV_p$ on $T_pM^{n+1}$.
Denote by $dA_p$ the induced volume form by $dV_p$ on
$S_pM^{n+1}$.

Define the diffeomorfism $\varphi_t:S_pM^{n+1}\rightarrow
S_t^{n}(p)$ [4] by
$$\varphi_t(y)=\exp_p(ty),\ y \in S_pM^{n+1}, \ t\geqslant0.$$

Let $dA_t$ denote the induced volume form on $S_t^{n}(p)$ by
$dV_F$. Define $$\eta_t:S_pM^{n+1}\rightarrow [0,\infty)$$ by
\begin{equation}\varphi_t^*dA_t|_{\varphi_t(y)}=\eta_t(y)dA_p|_y.\end{equation}

Integrating (4) over $S_pM^{n+1}$, we have
$$Area(S_t^{n}(p))=\int_{S_pM^{n+1}}\eta_t(y)dA_p.$$

Applying the co-area formula (1), we obtain
$$Vol(B_r^{n+1}(p))=\int_0^t\left(\int_{S_pM^{n+1}}\eta_s(y)dA_p\right)ds.$$

Remark that in the Riemannian case $\eta_t$ is the Jacobian of the
exponential map, and the explicit expression for the Jacobian
gives us all the necessary estimates. Unfortunately, the
integration of such estimates only leads to "coarse"$\ $estimates
for Finsler geometry.

Now, let us estimate $\eta_t$. For a small number $\varepsilon>0$
define the flow
\begin{equation}\phi_{\varepsilon}(x)=\varphi_{t+\varepsilon}\circ
\varphi^{-1}_t(x),\ x\in S_t^{n}(p).\end{equation}

For a point $x\in S_t^{n}(p)$, there is an open neighborhood $U$
of $x$ such that $\phi_{\varepsilon}$ is defined on $U$. The
Cartan-Hadamard theorem guarantees the non-existence of conjugate
points in all $M^{n+1}$, i.e., the existence of metric balls of
arbitrary radii.

As in the definition of mean curvature define
$\Theta(x,\varepsilon)$ by
$$\phi_{\varepsilon}^*dA_{s+\varepsilon}|_x=\Theta(x,\varepsilon)dA_s|_x.$$
Using (4), (5) we get
\begin{equation}\Theta(x,\varepsilon)=\frac{\eta_{t+\varepsilon}(y)}{\eta_t(y)},\
x=\varphi_t(y).\end{equation}

Let $\Pi_t$ denote the mean curvature or $S_t^{n}(p)$ at $x$ with
respect to the outward-pointing normal vector. From the definition
of mean curvature and (6), we have
\begin{equation}\Pi_{t}=\frac{\partial}{\partial \varepsilon}\left(\ln
\Theta(x,\varepsilon)
\right)|_{\varepsilon=0}=\frac{d}{dt}(\ln\eta_{t}(y)
).\end{equation}

Define $\chi_i(t)$ by
$$\chi_i(t)=\left(e^{-\delta_it}\frac{\sinh(k_it)}{k_i}\right)^n.$$

Then we have \begin{equation}
\frac{d}{dt}(\ln\chi_i(t))=nk_i\coth(k_it)-n\delta_i.\end{equation}

Taking into account the restrictions on curvature we can apply
Theorem 3. Then using (2), (3) we get
$$nk_1\coth(k_1t)-n\delta_1\leqslant \Pi_t\leqslant
nk_2\coth(k_2t)-n\delta_2.$$ This implies:
$$\frac{d}{dt}\left(\frac{\eta_{t}(y)}{\chi_2(t)}\right)\leqslant0, \
\frac{d}{dt}\left(\frac{\eta_{t}(y)}{\chi_1(t)}\right)\geqslant0,$$
and
$$\eta_{t_2}(y)\chi_1(t_1)\geqslant \eta_{t_1}(y)\chi_1(t_2),$$
$$\eta_{t_2}(y)\chi_2(t_1)\leqslant \eta_{t_1}(y)\chi_2(t_2),\ 0<t_1\leqslant t_2.$$

Integrating over $S_pM^{n+1}$ with respect to $dA_p$, we obtain
$$Area(S_{t_2}^{n}(p))\chi_1(t_1)\geqslant Area(S_{t_1}^{n}(p))\chi_1(t_2),$$
$$Area(S_{t_2}^{n}(p))\chi_2(t_1)\leqslant Area(S_{t_1}^{n}(p))\chi_2(t_2),\ 0<t_1\leqslant t_2.$$

Integrating from $0$ to $t_2$ with respect to $t_1$, we obtain
$$Area(S_{t_2}^{n}(p))\int_0^{t_2}\chi_1(t)dt\geqslant Vol(B_{t_2}^{n+1}(p))\chi_1(t_2),$$
$$Area(S_{t_2}^{n}(p))\int_0^{t_2}\chi_2(t)dt\leqslant Vol(B_{t_2}^{n+1}(p))\chi_2(t_2),\ 0<t_2.$$
Hence, we get
$$\frac{\chi_1(r)}{\int_0^{r}\chi_1(t)dt} \leqslant \frac{Area(S_{r}^{n}(p))}{Vol(B_{r}^{n+1}(p))}
\leqslant \frac{\chi_2(r)}{\int_0^{r}\chi_2(t)dt},$$ or
$$\frac{\int_0^{r}\left(e^{-\delta_2t}\sinh(k_2t)\right)^ndt}{\left(e^{-\delta_2r}\sinh(k_2r)\right)^n} \leqslant \frac{Vol(B_{r}^{n+1}(p))}{Area(S_{r}^{n}(p))}
\leqslant
\frac{\int_0^{r}\left(e^{-\delta_1t}\sinh(k_1t)\right)^ndt}{\left(e^{-\delta_1r}\sinh(k_1r)\right)^n},
\ r>0$$

Let us estimate these integrals.

$$\frac{\int_0^{r}\left(e^{-\delta_1t}\sinh(k_1t)\right)^ndt}{\left(e^{-\delta_1r}\sinh(k_1r)\right)^n}=
\frac{1}{(e^{-\delta_1r})^n}\int_0^r\left( e^{-\delta_1t}
\frac{e^{k_1t}-e^{-k_1t}}{e^{k_1r}-e^{-k_1r}}\right)^n dt
\leqslant$$
$$\leqslant \frac{1}{(e^{-\delta_1r})^n}
\int_0^r\left(e^{-\delta_1t+k_1(t-r)}\right)^ndt=\frac{e^{n\delta_1r}}{n(k_1-\delta_1)}\left(e^{-n\delta_1r}-e^{-nk_1r}\right)=$$
$$=\frac{1}{n(k_1-\delta_1)}\left(1-e^{-nr(k_1-\delta_1)}\right):=\mathcal{F}(r)$$

We can estimate the following integral by using the fact that
$(1-a)^n\geqslant 1-na$ for $0\leqslant a\leqslant 1$.

$$\frac{\int_0^{r}\left(e^{-\delta_2t}\sinh(k_2t)\right)^ndt}{\left(e^{-\delta_2r}\sinh(k_2r)\right)^n}=
\frac{e^{n\delta_2r}}{\left(1-e^{-2k_2r}\right)^n}\int_0^re^{-n\delta_2t}\left(1-e^{-2k_2t}\right)^ne^{k_2n(t-r)}dt\geqslant$$
$$\geqslant
\frac{e^{n\delta_2r}}{\left(1-e^{-2k_2r}\right)^n}\int_0^re^{-n\delta_2t}\left(1-ne^{-2k_2t}\right)e^{k_2n(t-r)}dt=$$
$$\frac{e^{n\delta_2r}}{\left(1-e^{-2k_2r}\right)^n}\left[\frac{1}{n(k_2-\delta_2)}\left(e^{-n\delta_2r}-e^{-nk_2r}\right)-\frac{n}{n(k_2-\delta_2)-2k_2}\left(e^{-n\delta_2r-2k_2r}-e^{-nk_2r}\right) \right]=$$
$$\frac{1}{\left(1-e^{-2k_2r}\right)^n}\left[\frac{1}{n(k_2-\delta_2)}\left(1-e^{-n(k_2-\delta_2)r}\right)-\frac{n}{(k_2-\delta_2)-2k_2}\left(e^{-2k_2r}-e^{-n(k_2-\delta_2)r}\right) \right]:=f(r) $$

Thus, we have
$$f(r) \leqslant \frac{Vol( B_r^{n+1}(p))}{Area(S_r^{n}(p))}\leqslant \mathcal{F}(r).$$

Using the inequalities $\delta_i<k_i$, we have

$$\lim_{r\rightarrow\infty}f(r)=\frac{1}{n(k_2-\delta_2)}$$

$$\lim_{r\rightarrow\infty}\mathcal{F}(r)=\frac{1}{n(k_1-\delta_1)}$$

As a consequence, we have

$$\frac{1}{n(k_2-\delta_2)}\leqslant \lim_{r\rightarrow\infty}\inf\frac{Vol( B_r^{n+1}(p))}{Area(S_r^{n}(p))}\leqslant
\lim_{r\rightarrow\infty}\sup\frac{Vol(
B_r^{n+1}(p))}{Area(S_r^{n}(p))}\leqslant
\frac{1}{n(k_1-\delta_1)}.$$

In the case when $K=-k^2$, $k>0$, $S=n\delta$, $\delta<k$, by
denoting $k_1=k_2=k$, $\delta_1=\delta_2=\delta,$ we have
$$\lim_{r\rightarrow\infty}\frac{Vol(B_r^{n+1}(p))}{Area(S_r^{n}(p))}=\frac{1}{n(k-\delta)}.$$

This completes the proof.

 $\blacksquare$

E x a m p l e 1. Let $U$ be a open bounded strongly convex domain
in $\mathbb{R}^n$. Take a point $x \in U$ and a direction $y \in
T_xU\backslash\{0\}\simeq U\backslash\{0\}$. Then the \textit{Funk
metric $F(x,y)$} is a Finsler metric that satisfies the following
condition
$$x+\frac{y}{F(x,y)}\in \partial U.$$
The indicatrix at each point for the Funk metric is a domain that
is a translate of $U$.

The \textit{Hilbert metric} is a symmetrized Funk metric: $$
\tilde{F}(x,y):=\frac{1}{2} \left (F(x,y)+F(x,-y) \right).$$

Note that for the Funk metric $B_x^n=U$. Thus
$$\sigma_F(x)=\frac{Vol_E(\mathbb{B}^n)}{Vol_E(B_x^n)}=\frac{Vol_E(\mathbb{B}^n)}{Vol_E(U)}=const.$$

Let $F$ be the Funk metric and let $\overline{F}$ be the Hilbert
metric on a strongly convex domain $U$ in $\mathbb{R}^n$.

Then geodesics of the Funk and Hilbert metrics are straight lines,
the Funk metric is of constant flag curvature $-\frac{1}{4}$, the
Hilbert metric is of constant flag curvature $-1$, and  the Funk
metric is of constant S-curvature $\frac{n+1}{2}$ [4].

Let $F$ be the Funk metric on a strongly convex domain $U$ in
$\mathbb{R}^{n+1}$. It is known that the $S$-curvature is equal to
$S=\frac{n+2}{2}=n\delta$, flag curvatures is equal to
$-k^2=-\frac{1}{4}$. Then the condition $\delta<k$ does not hold.

It is known that for the Funk metric
$$\frac{Vol(B_{r}^{n+1}(p))}{Area(S_{r}^{n}(p))}
=
\frac{\int_0^{r}\left(e^{-\frac{n+2}{2n}t}\sinh(\frac{t}{2})\right)^ndt}{\left(e^{-\frac{n+2}{2n}r}\sinh(\frac{r}{2})\right)^n}
$$
and one can show that
$$\lim_{r\rightarrow\infty}\frac{Vol(
B_r^{n+1}(p))}{Area(S_r^{n}(p))}=\infty.$$

Indeed, using Mathematica program, one can compute that
$$\frac{\int_0^{r}\left(e^{-\frac{n+2}{2n}t}\sinh(\frac{t}{2})\right)^ndt}{\left(e^{-\frac{n+2}{2n}r}\sinh(\frac{r}{2})\right)^n}=\frac{(e^r-1)}{n+1}\left(\frac{e^{-\frac{(n+1)}{n}r}(e^r-1)}{e^{-\frac{(n+r)}{n}r}(e^r-1)}\right)^n
$$

It is clear that such function grows to infinity as $r$ tends to
infinity.

In an $(n+1)$-dimensional Euclidean space such ratio also tends to
infinity.

This shows that the restrictions $\delta_i<k_i$ in the hypothesis
of the theorem are essential.

 $\square$

\section{Estimates on the volume growth entropy}

Let $(M^{n+1},F)$ be a Finsler manifold. Then the exponential
speed of the volume growth of a ball of radius $t>0$ is called the
\textit{volume growth entropy} of $(M^{n+1},F)$. The explicit
expression for the volume growth entropy is given by
$$\lim_{t\rightarrow \infty}\frac{\ln(Vol(B_{t}^{n+1}(p))}{t}.$$

In this section we estimate the volume growth entropy of a
Finsler-Hadamard manifold with the pinched flag curvature and the
$S$-curvature.

\textbf{Theorem 5.} \textit{Let $(M^{n+1},F)$ be an
$(n+1)$-dimensional Finsler-Hadamard manifold that satisfies the
following condition:
\begin{enumerate}
\item Flag curvature satisfies the inequalities $-k_2^2\leqslant K
\leqslant -k_1^2$, $k_1,k_2>0$, \item $S$-curvature satisfies the
inequalities $n\delta_1\leqslant S \leqslant n\delta_2$ such that
$\delta_i<k_i.$
\end{enumerate}}

\textit{Then we have}

$$n(k_1-\delta_1)\leqslant\lim_{t\rightarrow \infty}\frac{\ln(Vol(B_{t}^{n+1}(p))}{t}\leqslant n(k_2-\delta_2)$$

\textit{If $(M^{n+1},F)$ is a space of constant flag curvature
$K=-k^2$ and S-curvature $S=n\delta$, $\delta<k$, we have}
$$\lim_{t\rightarrow \infty}\frac{\ln(Vol(B_{t}^{n+1}(p))}{t}=n(k-\delta)$$

P r o o f $\ $ o f $\ $ T h e o r e m 5 :

Define $\chi_i(t)$ by
$$\chi_i(t)=\left(e^{-\delta_it}\frac{\sinh(k_it)}{k_i}\right)^n.$$

It was proved in [3], [4] that under the conditions 1. and 2. the
volume of a metric ball satisfies
\begin{equation}
Vol_E(\mathbb{S}^n)\int_0^t\chi_1(s)ds \leqslant
Vol(B_{t}^{n+1}(p))\leqslant
Vol_E(\mathbb{S}^n)\int_0^s\chi_2(s)ds
\end{equation}

By the direct computation, we have

$$\int_0^t\left(e^{-\delta_2s}\frac{\sinh(k_2s)}{k_2}\right)^n ds\leqslant\frac{1}{k_2^n}\int_0^t e^{sn(k_2-\delta_2)}ds=$$
$$=\frac{1}{n(k_2-\delta_2)k_2^n}\left(e^{tn(k_2-\delta_2)}-1\right)$$
Therefore, we get
$$\lim_{t\rightarrow \infty}\frac{\ln(Vol(B_{t}^{n+1}(p))}{t}\leqslant n(k_2-\delta_2).$$

Next,
$$\int_0^t\left(e^{-\delta_1s}\frac{\sinh(k_1s)}{k_1}\right)^n ds\geqslant\frac{1}{k_1^n}\int_0^t e^{-sn\delta_1}(1-ne^{-2k_1s})e^{k_1sn}ds=$$
$$=\frac{1}{k_1^n}\left[\frac{1}{n(k_1-\delta_1)}(e^{tn(k_1-\delta_1)}-1)+\frac{n}{k_1(n-2)-n\delta_1}(e^{tk_1(n-2)-n\delta_1}-1) \right].$$

This implies
$$\lim_{r\rightarrow \infty}\frac{\ln(Vol(B_{r}^{n+1}(p))}{r}\geqslant n(k_1-\delta_1)$$

And Theorem 5 follows easily.

$\blacksquare$

E x a m p l e 2. Let $F$ be the Funk metric on a strongly convex
domain $U$ in $\mathbb{R}^{n+1}$. Then the condition $\delta<k$
does not hold.

Then analogously as in Example 1 one can show that
$$\lim_{t\rightarrow \infty}\frac{\ln(Vol(B_{t}^{n+1}(p))}{t}=0.$$

In an (n+1)-dimensional Euclidean space such ratio also tends to
zero.

This shows that the restrictions $\delta_i<k_i$ in the hypothesis
of the theorem are essential.

In was shown in [7] that for the Hilbert metric $F$ on a strongly
convex domain $U$ in $\mathbb{R}^{n+1}$
$$\lim_{t\rightarrow \infty}\frac{\ln(Vol(B_{t}^{n+1}(p))}{t}=n.$$

Recall that $n$ is precisely the volume growth entropy of
$\mathbb{H}^{n+1}$.

 $\square$


\begin{thebibliography}{99}
\bibitem{Y1}
\textit{A.A. Borisenko, E. Gallego, A. Reventos.} Relation between
area and volume for convex sets in Hadamard manifolds. -
\textit{Differential geometry and its application} \textbf{14},
2001 - p. 267-280.
\bibitem{Y2}
\textit{A.A. Borisenko, V. Miquel.} Comparison theorem on convex
hypersurfaces in Hadamard manifolds. - \textit{Annals of Global
Analysis and Geometry} \textbf{21}, 2002 - p.191-202.
\bibitem{Y3}
\textit{Z. Shen.} Volume Comparison and its application in
Riemann-Finsler geometry. - \textit{Advanced in Math.}
\textbf{128}, 1997 - p.306-328..
\bibitem{Y4}
\textit{Z. Shen.} Lectures on Finsler Geometry. - \textit{World
Scientific Publishing Co}, 2001. - 306p.
\bibitem{Y5}
\textit{Z. Shen.} Landsberg curvature, S-curvature and Riemann
curvature. - \textit{Cambridge University Press}, 2004. - 53p.
\bibitem{Y6}
\textit{D. Egloff.} Uniform Finsler Hadamard manifolds. -
\textit{Annales de l'Institut Henri Poincare}, \textbf{Vol. 66},
\textbf{N 3}, 1997 - p. 323-357.
\bibitem{Y7}
\textit{B. Colbois, P. Verovic}. Hilbert Geometries for Stricly
Convex Domains. - \textit{Geometriae Dedicata} \textbf{105},  2004
- p. 29-42.
\end{thebibliography}
\end{document}